\documentclass[a4paper,11pt]{article}
\usepackage{amssymb}
\usepackage{amsfonts}
\usepackage{amsmath}

\input{tcilatex}

\begin{document}

\title{A K--theory approach to the tangent invariants of blow--ups}
\author{Haibao Duan\thanks{%
The author's research is supported by 973 Program 2011CB302400 and NSFC
11131008.} \\
%EndAName
Academy of Mathematics and Systems Sciences, \\
Chinese Academy of Sciences, Beijing 100190\\
dhb@math.ac.cn}
\date{}
\maketitle

\begin{abstract}
We extend the formula for the Chern classes of blow-ups of algebraic
varieties due to Porteous and Lascu--Scott, and of symplectic and complex
manifolds due to Geiges and Pasquotto, to the blow--ups of almost complex
manifolds.

Our approach is based on a concrete partition for the tangent bundle of a
blow--up. The use of topological K--theory of vector bundles simplifies the
classical approaches.

\begin{description}
\item[2000 Mathematical Subject Classification:] 53D35, 57R20

\item[Key words and phrases:] Blow ups; Almost complex manifolds; K--theory,
Chern classes.
\end{description}
\end{abstract}

\section{Introduction}

In this paper all manifolds under consideration are in the real and smooth
category, which are connected but not necessarily compact and orientable.

Let $X\subset M$\ be a smooth submanifold whose normal bundle has a complex
structure, and let $\widetilde{M}$ be the blow--up of $M$ along $X$. We
present a partition for the tangent bundle $\tau _{\widetilde{M}}$ of $%
\widetilde{M}$ which implies that, if $X\subset M$ is an embedding in the
category of \textsl{almost complex manifolds}, then the blow-up $\widetilde{M%
}$ has a canonical almost complex structure, see Theorem 2.3 and Corollary
2.4.

The partition on $\tau _{\widetilde{M}}$ is ready to apply to deduce a
formula for the tangent bundle $\tau _{\widetilde{M}}$ of the blow--up $%
\widetilde{M}$ in the $K$--theory $K(\widetilde{M})$ of complex bundles over 
$\widetilde{M}$, see Theorem 3.2 and Remark 3.3, which in turn yields a
formula for the total Chern class $C(\widetilde{M})$ of $\widetilde{M}$, see
Theorem 4.6.

Historically, the formula for the Chern class of a blow--up of a nonsingular
variety was conjectured by J. A. Todd \cite{[T]} and B. Segre \cite{[S]},
confirmed by I. R. Porteous \cite{[P]} and by Lascu--Scott \cite{[LS1],
[LS2]}, generalized to the blow ups of possibly singular varieties along
regularly embedded centers by Aluffi \cite{[A]}. It has also been extended
to the blow--ups in the categories of symplectic and complex manifolds by H.
Geiges and F. Pasquotto \cite{[GP]}. While establishing the formula in its
natural generality we demonstrate also an approach with deserved simplicity.

The author would like to thank P. Aluffi for bring his attention to the
topic.

\section{Geometry of a blow--up}

For an Euclidean vector bundle $\xi $ over a topological space $Y$ write $%
D(\xi )$ and $S(\xi )$ for the unit disk bundle and sphere bundle of $\xi $,
respectively. If $A\subset Y$ is a subspace (resp. if $f:W\rightarrow Y$ is
a map) write $\xi \mid A$ (resp. $f^{\ast }\xi $) for the restriction of $%
\xi $ to $A$ (resp. the induced bundle over $W$).

Assume throughout that $i_{X}:X\rightarrow M$ is a smooth embedding that
embeds $X$ as a closed subset of $M$, and whose normal bundle $\gamma _{X}$
has a fixed complex structure $J$. Furnish $M$ with an Riemannian metric so
that the induced metric on $\gamma _{X}$ is Hermitian in the sense of Milnor 
\cite[p.156]{[MS]}.

Let $\pi :E=\mathbb{P}(\gamma _{X})$ $\rightarrow X$ be the complex
projective bundle associated with $\gamma _{X}$. Since the tautological line
bundle $\lambda _{E}$ on $E$ is a subbundle of the induced bundle $\pi
^{\ast }\gamma _{X}$ we can formulate the composition

\begin{quote}
$G:D(\lambda _{E})\subset D(\pi ^{\ast }\gamma _{X})\overset{\widehat{\pi }}{%
\rightarrow }D(\gamma _{X})$,
\end{quote}

\noindent where $\widehat{\pi }$ is the obvious bundle map over $\pi $.
Regard $E\subset D(\lambda _{E})$\ and $X\subset D(\gamma _{X})$\ as the
zero sections of the corresponding disk bundles, respectively.

\bigskip

\noindent \textbf{Lemma 2.1.} \textsl{The map }$G$\textsl{\ agrees with the
projection }$\pi $\textsl{\ on }$E$\textsl{,} \textsl{and} \textsl{restricts
to a diffeomorphism }$D(\lambda _{E})\smallsetminus E\rightarrow D(\gamma
_{X})\smallsetminus X$\textsl{.}

\noindent \textbf{Proof.} Explicitly we have

\begin{quote}
$D(\lambda _{E})=\{(l,v)\in E\times \pi ^{\ast }\gamma _{X}\mid v\in l\in
E,\left\Vert v\right\Vert ^{2}\leq 1\}$,

$D(\gamma _{X})=\{(x,v)\in X\times \gamma _{X}\mid v\in \gamma _{X}\mid
x,\left\Vert v\right\Vert ^{2}\leq 1\}$.
\end{quote}

\noindent The inverse of $G$ on $D(\gamma _{X})\smallsetminus X$ is $%
(x,v)\rightarrow (\left\langle v\right\rangle ,v)$, where $v\in \gamma
_{X}\mid x$ with $v\neq 0$ and where $\left\langle v\right\rangle \in E$
denotes the complex line spanned by a non--zero vector $v$.$\square $

\bigskip

It follows from Lemma 2.1 that the map $G$ restricts to a diffeomorphism

\begin{quote}
$g=$ $G\mid S(\lambda _{E}):S(\lambda _{E})=\partial D(\lambda
_{E})\rightarrow S(\gamma _{X})=\partial D(\gamma _{X})$.
\end{quote}

\noindent With respect to the metric on $M$ identifying $D(\gamma _{X})$
with a tubular neighborhood of $X$ in $M$ we can formulate the adjoint
manifold

\begin{enumerate}
\item[(2.1)] $\widetilde{M}=(M\setminus \overset{\circ }{D(\gamma _{X})}%
)\cup _{g}D(\lambda _{E})$
\end{enumerate}

\noindent by gluing $D(\lambda _{E})$ to $(M\setminus \overset{\circ }{%
D(\gamma _{X})})$ along $S(\lambda _{E})$ using $g$. Moreover, piecing
together the identity on $M\setminus \overset{\circ }{D(\gamma _{X})}$ and
the map $G$ yields the smooth map

\begin{enumerate}
\item[(2.2)] $f:\widetilde{M}=(M\setminus \overset{\circ }{D(\gamma _{X})}%
)\cup _{g}D(\lambda _{E})\rightarrow M=(M\setminus \overset{\circ }{D(\gamma
_{X})})\cup _{id}D(\gamma _{X})$,
\end{enumerate}

\noindent which is known as the\textsl{\ blow--up} of $M$ along $X$ with 
\textsl{exceptional divisor} $E$ \cite{[M]}. Obvious but useful properties
of the map $f$ are listed below.

\bigskip

\noindent \textbf{Lemma 2.2.} \textsl{Let }$i_{E}:E\rightarrow \widetilde{M}$%
\textsl{\ (resp. }$i_{X}:X\rightarrow M$\textsl{) be the zero section of }$%
D(\lambda _{E})$\textsl{\ (resp. of }$D(\gamma _{X})$\textsl{)} \textsl{in
view of the decomposition (2.1). Then}

\begin{quote}
\textsl{i) the normal bundle of }$E$\textsl{\ in }$\widetilde{M}$ \textsl{is 
}$\lambda _{E}$\textsl{;}

\textsl{ii) }$f^{-1}(X)=E$\textsl{\ with }$f\circ i_{E}=$\textsl{\ }$%
i_{X}\circ \pi $\textsl{;}

\textsl{iii) }$f$\textsl{\ restricts to a diffeomorphism: }$\widetilde{M}%
\setminus E\rightarrow M\setminus X$\textsl{.}$\square $
\end{quote}

For a manifold $N$ write $\tau _{N}$ for its tangent bundle. One has then
the obvious bundle decompositions

\begin{center}
$\tau _{D(\lambda _{E})}\mid S(\lambda _{E})=\tau _{S(\lambda _{E})}\oplus 
\mathbb{R}(\alpha _{1})$; $\tau _{M\setminus \overset{\circ }{D(\gamma _{X})}%
}\mid S(\gamma _{X})=\tau _{S(\gamma _{X})}\oplus \mathbb{R}(\alpha _{2})$,
\end{center}

\noindent where $\alpha _{1}$ (resp. $\alpha _{2}$) is the outward (resp.
inward) unit normal field along the boundary $S(\lambda _{E})=\partial
D(\lambda _{E})$ (resp. $S(\gamma _{X})=\partial (M\setminus \overset{\circ }%
{D(\gamma _{X})})$) with $\mathbb{R}(\alpha _{i})$ the trivial real line
bundle spanned by the field $\alpha _{i}$. Moreover, if we let $\tau _{g}$
be the tangent map of the diffeomorphism $g$ then the decomposition (2.1) of 
$\widetilde{M}$ indicates the partition

\begin{enumerate}
\item[(2.3)] $\tau _{\widetilde{M}}=\tau _{M\setminus \overset{\circ }{%
D(\gamma _{X})}}\dbigcup\limits_{h}\tau _{D(\lambda _{E})}$,
\end{enumerate}

\noindent where the gluing diffeomorphism

\begin{quote}
$h:$ $\tau _{D(\lambda _{E})}\mid S(\lambda _{E})\rightarrow \tau
_{M\setminus \overset{\circ }{D(\gamma _{X})}}\mid S(\gamma _{X})$
\end{quote}

\noindent is the bundle map over $g$ with

\begin{quote}
$h(t\alpha _{1},u)=(t\alpha _{2},\tau _{g}(u))$, $u\in \tau _{S(\lambda
_{E})}$, $t\in \mathbb{R}$.
\end{quote}

\noindent Indeed, the two bundles $\tau _{D(\lambda _{E})}\mid S(\lambda
_{E})$ and $\tau _{M\setminus \overset{\circ }{D(\gamma _{X})}}\mid S(\gamma
_{X})$ admit more subtle decompositions with respect to them the gluing map $%
h$ in (2.3) admits a useful presentation.

Let $p_{E}:$\ $\lambda _{E}\rightarrow E$ and $p_{X}:$\ $\gamma
_{X}\rightarrow X$ be the obvious projections. The same notions will be
reserved for their restrictions to the subspaces $S(\lambda _{E})\subset
D(\lambda _{E})\subset \lambda _{E}$ and $S(\gamma _{X})\subset D(\gamma
_{X})\subset \gamma _{X}$, respectively.

For a topological space $Y$ write $1_{\mathbb{C}}$ (resp. $1_{\mathbb{R}}$)
for the trivial complex line bundle $Y\times \mathbb{C}$ (resp. the trivial
real line bundle $Y\times \mathbb{R}$) over $Y$. For a complex vector bundle 
$\xi $ write $\xi ^{r}$ for its real reduction. As example the
trivialization over\textsl{\ }$S(\lambda _{E})$

\begin{enumerate}
\item[(2.4)] $(p_{E}{}^{\ast }\lambda _{E}\mid S(\lambda _{E}))^{r}=\mathbb{R%
}(\alpha _{1})\oplus \mathbb{R}(J(\alpha _{1}))$
\end{enumerate}

\noindent indicates that $p_{E}{}^{\ast }\lambda _{E}\mid S(\lambda _{E})=1_{%
\mathbb{C}}$, where $J(\alpha _{1})$ is a unit tangent vector field on $%
S(\lambda _{E})$.

Let

\begin{quote}
$\widehat{g}:g^{\ast }(\tau _{M\setminus \overset{\circ }{D(\gamma _{X})}%
}\mid S(\gamma _{X}))\rightarrow \tau _{M\setminus \overset{\circ }{D(\gamma
_{X})}}\mid S(\gamma _{X})$
\end{quote}

\noindent be the induced bundle of $g$ over $S(\lambda _{E})$, and let

\begin{quote}
$\kappa :\tau _{D(\lambda _{E})}\mid S(\lambda _{E})\rightarrow g^{\ast
}(\tau _{M\setminus \overset{\circ }{D(\gamma _{X})}}\mid S(\gamma _{X}))$
\end{quote}

\noindent be the bundle isomorphism over the identity of $S(\lambda _{E})$
so that $h=\widehat{g}\circ \kappa $ \cite[ Lemma 3.1]{[MS]}. With respect
to the Hermitian metric induced from $\gamma _{X}$ one has the orthogonal
decomposition $\pi ^{\ast }\gamma _{X}=\lambda _{E}\oplus \lambda
_{E}^{\perp }$ in which $\lambda _{E}^{\perp }$\ denotes the orthogonal
complement of $\lambda _{E}$\ in $\pi ^{\ast }\gamma _{X}$.

\bigskip

\noindent \textbf{Theorem 2.3.} \textsl{The tangent bundle of the blow up }$%
\widetilde{M}$\textsl{\ has the partition}

\begin{quote}
$\tau _{\widetilde{M}}=\tau _{M\setminus \overset{\circ }{D(\gamma _{X})}%
}\dbigcup\limits_{\widehat{g}\circ \kappa }\tau _{D(\lambda _{E})}$,
\end{quote}

\noindent \textsl{in which}

\begin{quote}
\textsl{i) }$\tau _{D(\lambda _{E})}\mid S(\lambda _{E})=(\pi \circ
p_{E})^{\ast }\tau _{X}\oplus (p_{E}{}^{\ast }\lambda _{E})^{r}\oplus
p_{E}{}^{\ast }Hom(\lambda _{E},\lambda _{E}^{\perp })^{r}$\textsl{;}

\textsl{ii) }$g^{\ast }(\tau _{M\setminus \overset{\circ }{D(\gamma _{X})}%
}\mid S(\gamma _{X}))=(\pi \circ p_{E})^{\ast }\tau _{X}\oplus
(p_{E}{}^{\ast }\lambda _{E})^{r}\oplus p_{E}{}^{\ast }(\lambda _{E}^{\perp
})^{r}$\textsl{.}
\end{quote}

\noindent \textsl{Moreover, with respect to the decompositions i) and ii)
the bundle isomorphism }$\kappa $\textsl{\ is given by}

\begin{quote}
\textsl{a) }$\kappa \mid (\pi \circ p_{E})^{\ast }\tau _{X}=id;$

\textsl{b) }$\kappa \mid (p_{E}{}^{\ast }\lambda _{E})^{r}=id$\textsl{;}

\textsl{c) }$\kappa (b)=b(\alpha _{1})$\textsl{\ for }$b\in
Hom(p_{E}{}^{\ast }\lambda _{E},p_{E}{}^{\ast }\lambda _{E}^{\perp })^{r}$%
\textsl{.}
\end{quote}

\noindent \textbf{Proof. }It follows from the standard decompositions

\begin{quote}
$\tau _{E}=\pi ^{\ast }\tau _{X}\oplus Hom(\lambda _{E},\lambda _{E}^{\perp
})^{r}$, $\quad \tau _{D(\lambda _{E})}=(p_{E}{}^{\ast }\lambda
_{E})^{r}\oplus p_{E}{}^{\ast }\tau _{E}$
\end{quote}

\noindent that

\begin{enumerate}
\item[(2.5)] $\tau _{D(\lambda _{E})}=(p_{E}{}^{\ast }\lambda
_{E})^{r}\oplus (\pi \circ p_{E})^{\ast }\tau _{X}\oplus p_{E}{}^{\ast
}Hom(\lambda _{E},\lambda _{E}^{\perp })^{r}$.
\end{enumerate}

\noindent Similarly, it comes from

\begin{quote}
$\tau _{D(\gamma _{X})}=p_{X}{}^{\ast }\tau _{X}\oplus p_{X}{}^{\ast }\gamma
_{X}$, $\pi {}^{\ast }\gamma _{X}=\lambda _{E}\oplus \lambda _{E}^{\perp }$,
\end{quote}

\noindent as well as the definition of $f$ that

\begin{enumerate}
\item[(2.6)] $f^{\ast }\tau _{D(\gamma _{X})}=(p_{E}{}^{\ast }\lambda
_{E})^{r}\oplus (\pi \circ p_{E})^{\ast }\tau _{X}\oplus (p_{E}{}^{\ast
}\lambda _{E}^{\perp })^{r}$.
\end{enumerate}

\noindent One obtains the relations i) and ii) by restricting the
decomposition (2.5) and (2.6) to the subspace $S(\lambda _{E})\subset
D(\lambda _{E})$, respectively.

Finally, properties a), b), c) are transparent in view of the relation $h=%
\widehat{g}\circ \kappa $, together with the description of $g$ indicated in
the proof of Lemma 2.1.$\square $

\bigskip

A manifold $M$ is called \textsl{almost complex }if its tangent bundle is
furnished with a complex structure $J_{M}$ \cite[p.151]{[MS]}. Given two
almost complex manifolds $(X,J_{X})$ and $(M,J_{M})$ an embedding $i_{X}:$ $%
X\rightarrow M$ is called \textsl{almost complex }if $\tau _{X}$ is a
complex subbundle of the restricted bundle $\tau _{M}\mid X$. In this
situation the normal bundle $\gamma _{X}$ of $X$ has the canonical complex
structure $J$ induced from that on $\tau _{M}\mid X$ and that on $\tau _{X}$%
, hence the blow--up $\widetilde{M}$ of $M$ along $X$ is defined. Moreover,
in view of the decomposition (2.1) we note that

\begin{quote}
i) $J_{M}$ restricts to an almost complex structure on $M\setminus \overset{%
\circ }{D(\gamma _{X})}$;

ii) the tubular neighborhood $D(\lambda _{E})$ of $E$ in $\widetilde{M}$ has
the canonical almost complex structure so that as a complex vector bundle

$\qquad \tau _{D(\lambda _{E})}=(\pi \circ p_{E})^{\ast }\tau _{X}\oplus
p_{E}{}^{\ast }\lambda _{E}\oplus Hom(p_{E}{}^{\ast }\lambda
_{E},p_{E}{}^{\ast }\lambda _{E}^{\perp })$
\end{quote}

\noindent (compare this with (2.5)). Since with respect to the induced
complex structures on $\tau _{D(\lambda _{E})}\mid S(\lambda _{E})$ and on $%
\tau _{M\setminus \overset{\circ }{D(\gamma _{X})}}\mid S(\gamma _{X})$ the
clutching map $h$ in (2.3) is $\mathbb{C}$--linear by Theorem 2.3, one
obtains

\bigskip

\noindent \textbf{Corollary 2.4.} \textsl{If }$i_{X}:$\textsl{\ }$%
X\rightarrow M$\textsl{\ is an embedding in the category of almost complex
manifolds, then the blow--up }$\widetilde{M}$\textsl{\ has a canonical
almost complex structure that is compatible with that on }$M\setminus 
\overset{\circ }{D(\gamma _{X})}$\textsl{\ and that on }$D(\lambda _{E})$.$%
\square $

\bigskip

\noindent \textbf{Remark 2.5.} The analogue of Corollary 2.4 in the
symplectic setting is due to McDuff \cite[Section 3]{[M]}, which concludes
that if $i_{X}:$\textsl{\ }$X\rightarrow M$\textsl{\ }is an embedding of
symplectic manifolds, then the blow--up $\widetilde{M}$\ admits a symplectic
form which coincides with the one on $M\backslash X$ off the exceptional
divisor $E$.$\square $

\section{The tangent bundle of a blow--up}

For a topological space $Y$ let $K(Y)$ (resp. $\widetilde{K}(Y)$) be the $K$%
--theory (reduced $K$--theory) of complex vector bundles over $Y$. If $%
i_{X}: $\textsl{\ }$X\rightarrow M$ is an embedding in the category of
almost complex manifolds, then the blow--up $\widetilde{M}$ has a canonical
almost complex structure by Corollary 2.4. In particular, the difference $%
\tau _{\widetilde{M}}-f^{\ast }\tau _{M}$ can be regarded as an element of
the ring $\widetilde{K}(\widetilde{M})$. In Theorem 3.2 below we obtain a
formula expressing the element $\tau _{\widetilde{M}}-f^{\ast }\tau _{M}\in 
\widetilde{K}(\widetilde{M})$ in term of the decomposition $p_{E}^{\ast
}\gamma _{X}=p_{E}^{\ast }\lambda _{E}\oplus p_{E}^{\ast }\lambda
_{E}^{\perp }$.

For a relative CW--complex $(Y,A)$ the inclusion $j:(Y,\emptyset
)\rightarrow (Y;A)$ induces a homomorphism

\begin{enumerate}
\item[(2.1)] $j^{\ast }:K(Y;A)\rightarrow \widetilde{K}(Y)$,
\end{enumerate}

\noindent where $K(Y;A)$ is the relative $K$--group of the pair $(Y;A)$
defined by

\begin{enumerate}
\item[(3.2)] $K(Y;A)=:\widetilde{K}(Y/A)$.
\end{enumerate}

\noindent In addition, the group $K(Y;A)$ admits another description useful
in the subsequent calculation.

\bigskip

\noindent \textbf{Lemma 3.1} (\cite[Theorem 2.6.1]{[At]}). \textsl{Any
element in the group }$K(Y;A)$\textsl{\ can be represented by a triple }$%
[\xi ,\eta ;\alpha ]$\textsl{\ in which }$\xi $\textsl{\ and }$\eta $\textsl{%
\ are vector bundles over }$Y$\textsl{\ and }$\alpha :\xi \mid A\rightarrow $%
\textsl{\ }$\eta \mid A$\textsl{\ is a bundle isomorphism. }

\textsl{Moreover, with respect to this representation} \textsl{of the group} 
$K(Y;A)$ \textsl{one has}

\begin{quote}
\textsl{i) the triple }$[\xi ,\xi ;id]$\ \textsl{represents the zero} 
\textsl{for any bundle }$\xi $\textsl{\ over }$Y$\textsl{;}

\textsl{ii) }$[\xi ,\eta ;\alpha ]+[\xi _{1},\eta _{1};\alpha _{1}]=[\xi
\oplus \xi _{1},\eta \oplus \xi _{1};\alpha \oplus \alpha _{1}]$\textsl{;}

\textsl{iii) }$[\xi ,\eta ;\alpha ]\otimes \gamma =[\xi \otimes \gamma ,\eta
\otimes \gamma ;\alpha \otimes id]$\textsl{;}

\textsl{iv) }$j^{\ast }[\xi ,\eta ;\alpha ]=\xi -\eta $\textsl{,}
\end{quote}

\noindent \textsl{where }$\oplus $\textsl{\ means direct sum of vector
bundles (homomorphisms), and where }$\otimes $ \textsl{denotes the action} $%
K(X;A)\otimes K(X)\rightarrow K(X;A)$\textsl{\ defined by the tensor product
of vector bundles}.$\square $

\bigskip

For an embedding $i_{X}:X\rightarrow M$ of almost complex manifolds let $f:%
\widetilde{M}\rightarrow M$\ be the blow--up of $M$\ along $X$ with
exceptional divisor $E$. Consider the composition

\begin{quote}
$j_{E}:K(D(\lambda _{E}),S(\lambda _{E}))\underset{\cong }{\rightarrow }K(%
\widetilde{M},\widetilde{M}\smallsetminus \overset{\circ }{D}(\lambda _{E}))%
\overset{j^{\ast }}{\rightarrow }\widetilde{K}(\widetilde{M})$
\end{quote}

\noindent in which the first map is the \textsl{excision isomorphism}. By
Lemma 3.1 the trivialization $\varepsilon :\overline{\lambda }_{E}\mid
S(\lambda _{E})\rightarrow 1_{\mathbb{C}}$ indicated by (2.4) defines an
element $[p_{E}^{\ast }\overline{\lambda }_{E},1_{\mathbb{C}};\varepsilon
]\in K(D(\lambda _{E}),S(\lambda _{E}))$.

\bigskip

\noindent \textbf{Theorem 3.2.} \textsl{In the ring\ }$\widetilde{K}(%
\widetilde{M})$\textsl{\ one has}

\begin{enumerate}
\item[(3.3)] $\tau _{\widetilde{M}}-f^{\ast }\tau _{M}=j_{E}([p_{E}^{\ast }%
\overline{\lambda }_{E},1_{\mathbb{C}};\varepsilon ]\otimes p_{E}^{\ast
}\lambda _{E}^{\perp })$.
\end{enumerate}

\noindent \textbf{Proof.} The partition (2.1) of the blow up $\widetilde{M}$
implies the relation

\begin{quote}
$\tau _{\widetilde{M}}\mid M\backslash \overset{\circ }{D}(\gamma
_{X})=f^{\ast }\tau _{M}\mid M\backslash \overset{\circ }{D}(\gamma _{X})$
\end{quote}

\noindent which gives rise to an element $[\tau _{\widetilde{M}},f^{\ast
}\tau _{M};id]\in K(\widetilde{M};\widetilde{M}\smallsetminus \overset{\circ 
}{D(\lambda _{E})})$ by Lemma 3.1. Moreover, with respect to the \textsl{%
excision isomorphism}

\begin{quote}
$K(\widetilde{M};\widetilde{M}\smallsetminus \overset{\circ }{D(\lambda _{E})%
})\overset{\cong }{\rightarrow }K(D(\lambda _{E});S(\lambda _{E}))$
\end{quote}

\noindent we have $[\tau _{\widetilde{M}},f^{\ast }\tau _{M};id]=[\tau
_{D(\lambda _{E})},f^{\ast }\tau _{D(\gamma _{_{X}})};\kappa ]$, where

\begin{quote}
$\kappa :\tau _{D(\lambda _{E})}\mid S(\lambda _{E})\rightarrow g^{\ast
}(\tau _{M\setminus \overset{\circ }{D(\gamma _{X})}}\mid S(\gamma _{X}))$
\end{quote}

\noindent is the bundle isomorphism specified in Theorem 2.3. Granted with
the decompositions of the bundles $\tau _{D(\lambda _{E})}$ and $f^{\ast
}\tau _{D(\gamma _{X})}$ in (2.4) and (2.5), as well as the decomposition of
bundle map $\kappa $ in Theorem 2.3, we calculate

\begin{quote}
$\qquad \lbrack \tau _{D(\lambda _{E})},f^{\ast }\tau _{D(\gamma
_{X})};\kappa ]$

$=[(\pi \circ p_{E})^{\ast }\tau _{X},(\pi \circ p_{E})^{\ast }\tau
_{X};id]+[p_{E}{}^{\ast }\lambda _{E},p_{E}{}^{\ast }\lambda _{E};id]$

$\qquad +[p_{E}{}^{\ast }Hom(\lambda _{E},\lambda _{E}^{\perp
}),p_{E}{}^{\ast }\lambda _{E}^{\perp };\kappa ^{\prime }]$ (by ii) of Lemma
3.1)

$=[p_{E}{}^{\ast }Hom(\lambda _{E},\lambda _{E}^{\perp }),p_{E}{}^{\ast
}\lambda _{E}^{\perp };\kappa ^{\prime }]$ (by i) of Lemma 3.1)

$=[p_{E}{}^{\ast }(\overline{\lambda }_{E}\otimes \lambda _{E}^{\perp
}),p_{E}{}^{\ast }\lambda _{E}^{\perp };\kappa ^{\prime }]$ (since $%
Hom(\lambda _{E},\lambda _{E}^{\perp })=\overline{\lambda }_{E}\otimes
\lambda _{E}^{\perp }$)

$=[p_{E}{}^{\ast }\overline{\lambda }_{E},1_{\mathbb{C}};\varepsilon
]\otimes p_{E}{}^{\ast }\lambda _{E}^{\perp }$ (by iii) of Lemma 3.1)
\end{quote}

\noindent where $\kappa ^{\prime }$ is the restriction of $\kappa $ to the
direct summand $p_{E}{}^{\ast }Hom(\lambda _{E},\lambda _{E}^{\perp })$ of $%
\tau _{D(\lambda _{E})}$ (see c) of Theorem 2.3). Summarizing, in the group $%
K(D(\lambda _{E}),S(\lambda _{E}))$ we we have the relation

\begin{quote}
$[\tau _{\widetilde{M}},f^{\ast }\tau _{M};id]=[p_{E}{}^{\ast }\overline{%
\lambda }_{E},1_{\mathbb{C}};\varepsilon ]\otimes p_{E}{}^{\ast }\lambda
_{E}^{\perp }$.
\end{quote}

\noindent Finally, applying the map $j_{E}$ to both sides yields the formula
(3.3) by iv) of Lemma 3.1.$\square $

\bigskip

\noindent \textbf{Remark 3.3.} In the group $K(\widetilde{M})$ formula (3.3)
has the concise expression

\begin{quote}
$\tau _{\widetilde{M}}=f^{\ast }\tau _{M}+i_{E!}(\lambda _{E}^{\perp })$
\end{quote}

\noindent where $i_{E!}$ is the \textsl{Gysin map} in $K$--theory

\begin{quote}
$i_{E!}:K(E)\overset{\psi _{E}}{\underset{\cong }{\rightarrow }}K(D(\lambda
_{E}),S(\lambda _{E}))\overset{j_{E}}{\rightarrow }\widetilde{K}(\widetilde{M%
})$,
\end{quote}

\noindent and where $\psi _{E}$ is the \textsl{Thom isomorphism }$%
x\rightarrow $ $\psi _{E}(x)=U_{E}\otimes p_{E}^{\ast }x$, $x\in K(E)$, with 
$U_{E}\in K(D(\lambda _{E}),S(\lambda _{E}))$ the \textsl{Thom class} of $%
\lambda _{E}$. Indeed, from the general construction \cite[p.98--99]{[At]}
of Thom classes from the exterior algebras of vector bundles one finds that $%
U_{E}=[p_{E}^{\ast }\overline{\lambda }_{E},1;\varepsilon ]$.$\square $

\section{The Chern class of a blow up}

Based on the formula (3.3) we deduce a formula for the total Chern class $C(%
\widetilde{M})$ of a blow up $\widetilde{M}$ in the category of almost
complex manifolds. In this section the coefficients for cohomologies will be
the ring $\mathbb{Z}$ of integers.

\subsection{Preliminaries in Chern classes}

Let $BU$ be the classifying space of stable equivalent classes of complex
vector bundles, and let $c_{r}\in H^{2r}(BU)$ be the $r^{th}$ Chern class of
the universal complex vector bundle over $BU$. Then

\begin{quote}
$H^{\ast }(BU)=\mathbb{Z}[c_{1},c_{2},\cdots ]$.
\end{quote}

For a topological space $Y$ let $[Y,BU]$ be the set of homotopy classes of
maps from $Y$ to $BU$. In view of the standard identification

\begin{quote}
$\widetilde{K}(Y)=[Y,BU]$
\end{quote}

\noindent (\cite[p.210]{[Sw]}) we can introduce the total Chern class for
elements in $\widetilde{K}(Y)$ as the co--functor $C:\widetilde{K}%
(Y)\rightarrow H^{\ast }(Y)$ defined

\begin{quote}
$C(\beta )=1+f^{\ast }c_{1}+f^{\ast }c_{2}+\cdots $, $\beta \in \widetilde{K}%
(Y),$
\end{quote}

\noindent where $f:Y\rightarrow BU$ is the classifying map of the element $%
\beta $. Clearly one has

\bigskip

\noindent \textbf{Lemma 4.1.} \textsl{The transformation }$C$\textsl{\
satisfies the next two properties.}

\textsl{i) If }$\xi _{i}$\textsl{, }$i=1,2$\textsl{, are two complex vector
bundles over }$Y$\textsl{\ with equal dimension and with (the usual) total
Chern classes }$C(\xi _{i})$\textsl{,} \textsl{then}

\begin{quote}
$C(\xi _{1}-\xi _{2})=C(\xi _{1})C(\xi _{2})^{-1}$.
\end{quote}

\textsl{ii) For a closed subspace }$A\subset Y$\textsl{\ let }$%
j_{A}:(Y,\emptyset )\rightarrow (Y,A)$\textsl{\ and }$q_{A}:Y\rightarrow Y/A$%
\textsl{\ be the inclusion and quotient maps, respectively. Then the next
diagram commutes:}

\begin{quote}
$%
\begin{array}{ccc}
K(Y;A)=\widetilde{K}(Y/A) & \overset{j_{A}^{\ast }}{\rightarrow } & 
\widetilde{K}(Y) \\ 
C\downarrow &  & \downarrow C \\ 
H^{\ast }(Y/A) & \overset{q_{A}^{\ast }}{\rightarrow } & H^{\ast }(Y)%
\end{array}%
$.$\square $
\end{quote}

For two complex vector bundles $\lambda $ and $\xi $ over a space $Y$ with $%
\dim \lambda =1$, $\dim \xi =m$, and with the total Chern classes

\begin{quote}
$C(\lambda )=1+t$;$\quad C(\xi )=1+c_{1}(\xi )+c_{2}(\xi )+\cdots +c_{m}(\xi
)$,
\end{quote}

\noindent respectively, assume that the Chern roots of $\xi $ is $%
s_{1},\cdots ,s_{m}$. That is

\begin{quote}
$C(\xi )=\tprod_{1\leq i\leq m}(1+s_{i})$
\end{quote}

\noindent with $c_{r}(\xi )=e_{r}(s_{1},\cdots ,s_{m})$ the $r^{th}$
elementary symmetric function in the roots $s_{1},\cdots ,s_{m}$. The
calculation

\begin{quote}
$C(\lambda \otimes \xi )=\tprod\limits_{1\leq i\leq
m}(1+t+s_{i})=(1+t)^{m}\tprod\limits_{1\leq i\leq m}(1+\frac{s_{i}}{1+t})$

$\qquad =(1+t)^{m}[1+\frac{c_{1}(\xi )}{(1+t)}+\frac{c_{2}(\xi )}{(1+t)^{2}}%
+\cdots +\frac{c_{m}(\xi )}{(1+t)^{m}}]$
\end{quote}

\noindent shows that

\begin{enumerate}
\item[(4.1)] $C(\lambda \otimes \xi )=\dsum\limits_{0\leq r\leq
m}(1+t)^{m-r}c_{r}(\xi )$.
\end{enumerate}

For an Euclidean complex line bundle $\lambda $ over $Y$ with associated
disk bundle $p_{\lambda }:D(\lambda )\rightarrow Y$ the Thom space $%
T(\lambda )$ of $\lambda $ is the quotient space $D(\lambda )/S(\lambda )$.
In term of Lemma 3.1 the trivialization $\varepsilon :$ $p_{\lambda }^{\ast
}(\overline{\lambda })\mid S(\overline{\lambda })\rightarrow 1_{\mathbb{C}}$

\noindent indicated by (2.4) defines the element

\begin{quote}
$[p_{\lambda }^{\ast }(\overline{\lambda }),1_{\mathbb{C}};\varepsilon ]\in
K(D(\lambda ),S(\lambda ))=\widetilde{K}(T(\lambda ))$.
\end{quote}

\noindent Given a ring $A$ and a set $\{u_{1},\cdots ,u_{k}\}$ of $k$
elements let $A\{u_{1},\cdots ,u_{k}\}$ be the free $A$ module with basis $%
\{u_{1},\cdots ,u_{k}\}$.

In addtion to Theorem 3.2, the next result will play a key role in
establishing the formula for Chern class.

\bigskip

\noindent \textbf{Lemma 4.2. }\textsl{Let }$e\in H^{2}(Y)$\textsl{\ and }$%
x\in H^{2}(T(\lambda ))$\textsl{\ be respectively the Euler class and the
Thom class of the oriented bundle }$\lambda $\textsl{. Then }

\textsl{i) the integral cohomology ring of }$T(\lambda )$\textsl{\ is
determined by the additive presentation}

\begin{enumerate}
\item[(4.2)] $H^{\ast }(T(\lambda ))=\mathbb{Z}\oplus H^{\ast }(Y)\{x\}$%
\textsl{,\ }
\end{enumerate}

\noindent \textsl{together with the single relation }$x^{2}+xe=0$\textsl{;}

\textsl{ii) for a complex vector bundle }$\xi $ \textsl{over }$Y$\textsl{\
with total Chern class }$C(\xi )=1+c_{1}+\cdots +c_{m}\in H^{\ast }(Y),$ 
\textsl{the total Chern class of the element }$[p_{\lambda }^{\ast }(%
\overline{\lambda }),1_{\mathbb{C}};\varepsilon ]\otimes p_{\lambda }^{\ast
}\xi \in \widetilde{K}(T(\lambda ))$\textsl{\ is}

\begin{enumerate}
\item[(4.3)] $C([p_{\lambda }^{\ast }(\overline{\lambda }),1_{\mathbb{C}%
};\varepsilon ]\otimes p_{\lambda }^{\ast }\xi )=(\dsum\limits_{0\leq r\leq
m}(1+x)^{m-r}c_{r})C(\xi )^{-1}\in H^{\ast }(T(\lambda )).$
\end{enumerate}

\noindent \textbf{Proof.} The presentation (4.2) comes immediately from the
Thom isomorphism theorem, which states that product with Thom class $x$
yields an additive isomorphism

\begin{quote}
$H^{r}(Y)\cong H^{r+2}(T(\lambda ))$,$\quad y\rightarrow y\cdot x$,$\quad
y\in H^{r}(X)$
\end{quote}

\noindent for all $r\geq 0$. For i) it remains to justify the relation $%
x^{2}+xe=0$.

Let $p:S(\lambda \oplus 1_{\mathbb{R}})\rightarrow Y$ be the sphere bundle
of the Euclidean bundle $\lambda \oplus 1_{\mathbb{R}}$ and set

\begin{quote}
$D_{+(-)}(\lambda )=\{(u,t)\in S(\lambda \oplus 1_{\mathbb{R}})\mid t\geq 0$ 
$(t\leq 0)\}$.
\end{quote}

\noindent It is clear that

\begin{quote}
a) $S(\lambda \oplus 1_{\mathbb{R}})=D_{-}(\lambda )\cup D_{+}(\lambda )$, $%
S(\lambda )=D_{-}(\lambda )\cap D_{+}(\lambda )$,

b) both $D_{\pm }(\lambda )$ can be identified with the disk bundle $%
D(\lambda )$ of $\lambda $.
\end{quote}

\noindent In view of a) we have the canonical map onto the Thom space $%
T(\lambda )$

\begin{quote}
$q:S(\lambda \oplus 1_{\mathbb{R}})\rightarrow T(\lambda )=S(\lambda \oplus
1_{\mathbb{R}})/D_{-}(\lambda )$.
\end{quote}

\noindent Moreover, letting $u=q^{\ast }x\in H^{2}(S(\lambda \oplus 1_{%
\mathbb{R}}))$ the ring $H^{\ast }(S(\lambda \oplus 1_{\mathbb{R}}))$ has
the presentation (see \cite[Lemma 4]{[D]})

\begin{quote}
$H^{\ast }(S(\lambda \oplus 1_{\mathbb{R}}))=H^{\ast
}(Y)\{1,u\}/\left\langle u^{2}+ue\right\rangle $.
\end{quote}

\noindent The relation $x^{2}+xe=0$ on $H^{\ast }(T(\lambda ))$ is verified
by the relation $u^{2}+ue=0$ on $H^{\ast }(S(\lambda \oplus 1_{\mathbb{R}}))$%
, together with the fact that the induced ring map $q^{\ast }$ is
monomorphic onto the direct summand $\mathbb{Z}\oplus H^{\ast }(Y)\{u\}$ of $%
H^{\ast }(S(\lambda \oplus 1_{\mathbb{R}})$.

For ii) define over $S(\lambda \oplus 1_{\mathbb{R}})$ the complex line
bundle $\lambda _{u}$ by

\begin{quote}
$\lambda _{u}=p^{\ast }\overline{\lambda }\mid D_{+}(\lambda )\cup
_{\varepsilon }1_{\mathbb{C}}\mid D_{-}(\lambda )$ (in view of the partition
a)).
\end{quote}

\noindent Then $C(\lambda _{u})=1+u\in H^{\ast }(S(\lambda \oplus 1_{\mathbb{%
R}})$. Moreover, for a vector bundle $\xi $ over\textsl{\ }$Y$ the element

\begin{quote}
$[\lambda _{u},1_{\mathbb{C}};\varepsilon ]\otimes p^{\ast }\xi \in
K(S(\lambda \oplus 1_{\mathbb{R}}),D_{-}(\lambda ))$
\end{quote}

\noindent corresponds to the element $[(p_{\lambda }^{\ast }(\overline{%
\lambda }),1_{\mathbb{C}};\varepsilon )\otimes p_{\lambda }^{\ast }\xi ]\in
K(D(\lambda ),S(\lambda ))$ under the excision isomorphism $K(S(\lambda
\oplus 1_{\mathbb{R}}),D_{-}(\lambda ))\cong K(D(\lambda ),S(\lambda ))$
indicated by b), which is also mapped to the element

\begin{quote}
$\lambda _{u}\otimes p^{\ast }\xi -p^{\ast }\xi \in \widetilde{K}(S(\lambda
\oplus 1_{\mathbb{R}}))$
\end{quote}

\noindent under the induced homomorphism $j^{\ast }$ of the inclusion $%
j:(S(\lambda \oplus 1_{\mathbb{R}}),\emptyset )\rightarrow (S(\lambda \oplus
1_{\mathbb{R}}),D_{-}(\lambda ))$ by iv) of Lemma 3.1. It follows from i) of
Lemma 4.1 that

\begin{quote}
$C(j^{\ast }([(\lambda _{u},1_{\mathbb{C}};\varepsilon )\otimes p^{\ast }\xi
]))=C(\lambda _{u}\otimes p^{\ast }\xi )C(p^{\ast }\xi )^{-1}$

$=(\dsum\limits_{0\leq r\leq m}(1+u)^{m-r}c_{r})C(\xi )^{-1}$ (by the
formula (4.1)).
\end{quote}

\noindent The commutativity of the diagram in ii) of Lemma 4.1 implies that

\begin{quote}
$q^{\ast }C([(\lambda _{u},1_{\mathbb{C}};\varepsilon )\otimes p^{\ast }\xi
])=(\dsum\limits_{0\leq r\leq m}(1+u)^{m-r}c_{r})C(\xi )^{-1}$.
\end{quote}

\noindent One obtains the formula (4.3) from $q^{\ast }(x)=u$ and the
injectivity of $q^{\ast }$.$\square $

\subsection{The integral cohomology ring of a blow--up}

Let $f:\widetilde{M}\rightarrow M$ be the blow up of $M$ along a submanifold 
$i_{X}:X\rightarrow M$ whose normal bundle $\gamma _{X}$ has a complex
structure and with total Chern class $C(\gamma _{X})=1+c_{1}+\cdots
+c_{k}\in H^{\ast }(X),$ $k=\dim _{\mathbb{C}}\gamma _{X}$. Regard $%
D(\lambda _{E})$ as a normal disk bundle of the exceptional divisor $E$ and
consider the quotient map onto the Thom space of $\lambda _{E}$

\begin{quote}
$q:\widetilde{M}\rightarrow \widetilde{M}/(M\backslash \overset{\circ }{D}%
(\gamma _{X})=T(\lambda _{E})$.
\end{quote}

\noindent According to i) of Lemma 4.2 one has

\bigskip

\noindent \textbf{Lemma 4.3. }\textsl{Let }$t\in H^{2}(E)$\textsl{\ and }$%
x\in H^{2}(T(\lambda _{E}))$\textsl{\ be the Euler class and Thom class of
the oriented bundle }$\lambda _{E}$\textsl{, respectively. Then the ring }$%
H^{\ast }(T(\lambda _{E}))$ \textsl{is determined by the additive
presentation}

\begin{enumerate}
\item[(4.5)] $H^{\ast }(T(\lambda _{E}))=\mathbb{Z}\oplus H^{\ast
}(X)[x,tx,\cdots ,t^{k-1}x]$,
\end{enumerate}

\noindent \textsl{together with the relations}

\begin{quote}
\textsl{i) }$x^{2}+tx=0$\textsl{; ii) }$t^{k}+c_{1}\cdot t^{k-1}+\cdots
+c_{k-1}\cdot t+c_{k}=0$\textsl{.}$\square $
\end{quote}

Consider the Gysin maps of the embeddings $i_{X}:X\rightarrow M$ and $%
i_{E}:E\rightarrow \widetilde{M}$ in cohomology

\begin{center}
$i_{X!}:H^{\ast }(X)\overset{\psi _{X}}{\underset{\cong }{\rightarrow }}%
H^{\ast }(D(\gamma _{X}),S(\gamma _{X}))\underset{\cong }{\rightarrow }%
H^{\ast }(M,M\smallsetminus \overset{\circ }{D}(\gamma _{X}))\overset{%
j^{\ast }}{\rightarrow }H^{\ast }(M)$

$i_{E!}:H^{\ast }(E)\overset{\psi _{E}}{\underset{\cong }{\rightarrow }}%
H^{\ast }(D(\lambda _{E}),S(\lambda _{E}))\underset{\cong }{\rightarrow }%
H^{\ast }(\widetilde{M},\widetilde{M}\smallsetminus \overset{\circ }{D}%
(\lambda _{E}))\overset{j^{\ast }}{\rightarrow }H^{\ast }(\widetilde{M})$
\end{center}

\noindent where $\psi _{X}$ and $\psi _{E}$ are the Thom isomorphisms. We
shall set

\begin{quote}
$\omega _{E}=i_{E!}(1)\in H^{2}(\widetilde{M})$, $\omega _{X}=i_{X!}(1)\in
H^{2k}(M)$.
\end{quote}

\noindent Geometrically, if $M$ is closed and oriented the class $\omega
_{X} $ (for instance) is the Poincare dual of the cycle class $i_{X\ast
}[X]\in H_{\ast }(M)$.

The next result is shown in \cite[Theorem 1]{[DL]}.

\bigskip

\noindent \textbf{Theorem 4.4.} \textsl{The ring map }$f^{\ast }:H^{\ast
}(M)\rightarrow H^{\ast }(\widetilde{M})$\textsl{\ embeds the ring }$H^{\ast
}(M)$\textsl{\ as a direct summand of }$H^{\ast }(\widetilde{M})$,\textsl{\
and induces the decomposition}

\begin{enumerate}
\item[(4.6)] $H^{\ast }(\widetilde{M})=f^{\ast }(H^{\ast }(M))\oplus H^{\ast
}(X)\{\omega _{E},\cdots ,\omega _{E}^{k-1}\},$ $2k=\dim _{\mathbb{R}}\gamma
_{X}$
\end{enumerate}

\noindent \textsl{that is subject to the two relations}

\begin{quote}
\textsl{i) }$f^{\ast }(\omega _{X})=\tsum\limits_{1\leq r\leq
k}(-1)^{r-1}c_{k-r}\cdot \omega _{E}^{r}$;

\textsl{ii)} $f^{\ast }(y)\cdot \omega _{E}=i_{X}^{\ast }(y)\cdot \omega
_{E} $\textsl{, }$y\in H^{r}(M)$.
\end{quote}

\textsl{Moreover, with respect to the presentations of the rings }$H^{\ast
}(T(\lambda _{E}))$ \textsl{and} $H^{\ast }(\widetilde{M})$ \textsl{in (4.5)
and (4.6), the induced ring map }$q^{\ast }$\textsl{\ is determined by}

\begin{quote}
\textsl{iii)} $q^{\ast }(t^{r}x)=(-1)^{r}\omega _{E}^{r+1}$\textsl{,} $r\geq
0$. $\square $
\end{quote}

\noindent \textbf{Remark 4.5.} In \cite[p.605]{[GH]} Griffiths and Harris
obtained the decomposition (4.6) for blow ups of complex manifolds, while
the relations i) and ii) were absent. In the non--algebraic settings partial
information on the ring $H^{\ast }(\widetilde{M})$ was also obtained by
McDuff in \cite[Proposition 2.4]{[M]}.

In comparison the structure of $H^{\ast }(\widetilde{M})$ as a ring is
completely determined by the additive decomposition (4.6), together with the
relations i) and ii). Indeed, granted with the fact that $f^{\ast }(H^{\ast
}(M))\subset H^{\ast }(\widetilde{M})$ is a subring, the relations i) and
ii) are sufficient to express, respectively, the products of elements in the
second summand, and the products between elements in the first and second
summands, as elements in the decomposition (4.6). This idea has been applied
in \cite{[DL]} to determine the integral cohomology rings of the varieties
of complete conics and complete quadrics in the $3$--space $\mathbb{P}^{3}$,
and justify two enumerative problems due to Schubert.$\square $

\subsection{The Chern class of a blow up}

Assume now that $f:\widetilde{M}\rightarrow M$ is a blow up in the category
of almost complex manifolds, and that the total Chern class of the normal
bundle $\gamma _{X}$ is

\begin{quote}
$C(\gamma _{X})=1+c_{1}+\cdots +c_{k}$, $\dim _{\mathbb{R}}\gamma _{X}=2k$.
\end{quote}

\noindent Applying the transformation $C:\widetilde{K}(\widetilde{M}%
)\rightarrow H^{\ast }(\widetilde{M})$ to the equality (3.3) and noting the
obvious relation

\begin{quote}
$p_{E}^{\ast }\lambda _{E}^{\perp }=p_{E}^{\ast }\gamma _{X}-p_{E}^{\ast
}\lambda _{E}$
\end{quote}

\noindent in the $K$--group $K(D(\lambda _{E}))$ one obtains

\begin{quote}
$C(\widetilde{M})\cdot f^{\ast }C(M)^{-1}=C(j_{E}([p_{E}^{\ast }\overline{%
\lambda }_{E},1_{\mathbb{C}};\varepsilon ]\otimes (p_{E}^{\ast }\gamma
_{X}-p_{E}^{\ast }\lambda _{E}))$

$=\frac{C(j_{E}([p_{E}^{\ast }\overline{\lambda }_{E},1_{\mathbb{C}%
};\varepsilon ]\otimes p_{E}^{\ast }\gamma _{X}))}{C(j_{E}([p_{E}^{\ast }%
\overline{\lambda }_{E},1_{\mathbb{C}};\varepsilon ]\otimes p_{E}^{\ast
}\lambda _{E}))}$ (by i) of Lemma 4.1)

$=q^{\ast }\frac{C([p_{E}^{\ast }\overline{\lambda }_{E},1_{\mathbb{C}%
};\varepsilon ]\otimes p_{E}^{\ast }\gamma _{X})}{C([p_{E}^{\ast }\overline{%
\lambda }_{E},1_{\mathbb{C}};\varepsilon ]\otimes p_{E}^{\ast }\lambda _{E})}
$ (by ii) of Lemma 4.1).
\end{quote}

\noindent Furthermore, from

\begin{quote}
$C([p_{E}^{\ast }\overline{\lambda }_{E},1_{\mathbb{C}};\varepsilon ]\otimes
p_{E}^{\ast }\gamma _{X})=(\dsum\limits_{0\leq r\leq
k}(1+x)^{k-r}c_{r})C(\gamma _{X})^{-1}$,

$C([p_{E}^{\ast }\overline{\lambda }_{E},1_{\mathbb{C}};\varepsilon ]\otimes
p_{E}^{\ast }\lambda _{E})=(1+x+t)(1+t)^{-1}$
\end{quote}

\noindent by (4.3) one gets

\begin{quote}
$C(\widetilde{M})\cdot f^{\ast }C(M)^{-1}=q^{\ast }(\dsum\limits_{0\leq
r\leq k}(1+x)^{k-r}c_{r})(1+t)(1+x+t)^{-1}C(\gamma _{X})^{-1}$

$=(\dsum\limits_{0\leq r\leq k}(1+\omega _{E})^{k-r}c_{r})(1-\omega
_{E})C(\gamma _{X})^{-1}$.
\end{quote}

\noindent where the second equality comes from iii) of Theorem 4.4. It
implies that

\begin{quote}
$C(\widetilde{M})-f^{\ast }C(M)=f^{\ast }C(M)\cdot g(\omega _{E})$
\end{quote}

\noindent with

\begin{quote}
$g(\omega _{E})=(\dsum\limits_{0\leq r\leq k}(1+\omega
_{E})^{k-r}c_{r})(1-\omega _{E})C(\gamma _{X})^{-1}-1\in H^{\ast }(%
\widetilde{M})$
\end{quote}

\noindent a polynomial in $\omega _{E}$ with coefficients in $H^{\ast }(X)$.

It is crucial for us to observe from the obvious relation $g(0)=0$ that the
polynomial $g(\omega _{E})$ is divisible by $\omega _{E}$. Therefore, the
relation ii) in Theorem 4.4 is applicable to yield

\begin{quote}
$C(\widetilde{M})-f^{\ast }C(M)=f^{\ast }C(M)\cdot g(\omega _{E})$

$=i_{X}^{\ast }C(M)g(\omega _{E})$

$=C(X)C(\gamma _{X})g(\omega _{E})$ (since $\tau _{M}\mid X=\tau _{X}\oplus
\gamma _{X}$).
\end{quote}

\noindent Summarizing we get

\bigskip

\noindent \textbf{Theorem 4.6. }\textsl{With respect to decomposition of the
ring }$H^{\ast }(\widetilde{M})$\textsl{\ in (4.6), the total Chern class of
the blow up }$\widetilde{M}$\textsl{\ is }

\begin{quote}
$C(\widetilde{M})=f^{\ast }C(M)+C(X)(\dsum\limits_{0\leq r\leq k}(1+\omega
_{E})^{k-r}c_{r})(1-\omega _{E})-\dsum\limits_{0\leq r\leq k}c_{r})$.$%
\square $
\end{quote}

\noindent \textbf{Examples 4.7.} The simplest example assuring that Theorem
4.6 is non--trivial in the category of almost complex manifolds is the
following one. Recall that the $6$--dimensional sphere $S^{6}$ has a
canonical almost complex structure. The blow--up of $S^{6}$ at a point $X\in
S^{6}$ is diffeomorphic to the complex projective $3$--space $\mathbb{P}^{3}$%
, together with an induced almost complex structure $J$ (see Corollary 2.4).
By Theorem 4.6 we have

\begin{quote}
$C(\mathbb{P}^{3},J)=1-2x-4x^{3}$, 
\end{quote}

\noindent where $x\in H^{2}(\mathbb{P}^{3})$ is the exceptional divisor.
This computation shows that $J$ is different with the canonical complex
structure on $\mathbb{P}^{3}$.

Similarly, in the recent work \cite{[Y]} H. Yang classified those $(2n-1)$%
--connected $4n$--manifolds that admit almost complex structures. In term of
the Wall's invariant $\alpha :H_{2n}(M)\rightarrow \pi _{2n-1}(SO(2n))$ for
such a manifold $M$ \cite{[W]} the set of isotopy classes of smooth
embeddings $S^{2n}\rightarrow M$ whose normal bundles have complex
structures can be identified with an explicit subset $J(M)$ of the homology
group $H_{2n}(M)$, where $S^{2n}$ is the $2n$--dimensional sphere. Blow--ups
of $M$ along elements in the set $J(M)$ can provide us with a rich family of
almost complex manifolds of dimension $4n$, to which Therom 4.6 is directly
applicable to compute their total Chern classes.

\bigskip

\end{document}